\documentclass[11pt,reqno]{amsart}

\usepackage[latin1]{inputenc}
\usepackage[left=3.5cm,right=3.5cm,top=3.5cm,bottom=3.5cm]{geometry}
\usepackage{amssymb}
\usepackage{graphicx}
\usepackage{amscd}
\usepackage[hidelinks]{hyperref}
\usepackage{color}
\usepackage{float}
\usepackage{graphics,amsmath,amssymb}
\usepackage{amsthm}
\usepackage{amsfonts}
\usepackage{latexsym}
\usepackage{epsf}
\usepackage{enumerate}
\usepackage{xifthen}
\usepackage{mathrsfs}
\usepackage{dsfont}
\usepackage{makecell}
\usepackage[FIGTOPCAP]{subfigure}
\usepackage{amsmath}
\allowdisplaybreaks[4]
\usepackage{listings}
\usepackage{etoolbox}

 \newtheoremstyle{mytheorem}
 {3pt}
 {3pt}
 {\slshape}
 {}
 {\bfseries}
 {.}
 { }
 {}

\numberwithin{equation}{section}

\theoremstyle{mytheorem}
\newtheorem{theorem}{Theorem}[section]

\newtheorem{lemma}[theorem]{Lemma}

\theoremstyle{definition}

\newcommand{\doi}[1]{\href{http://dx.doi.org/#1}{doi: #1}}

\newcommand{\Keywords}[1]{\ifthenelse{\isempty{#1}}{}{\smallskip \smallskip \noindent \textbf{Keywords}. #1}}
\newcommand{\MSC}[2][2010]{\ifthenelse{\isempty{#2}}{}{\smallskip \smallskip \noindent \textbf{#1MSC}. #2}}

\makeatletter
\patchcmd{\@settitle}{\uppercasenonmath\@title}{\LARGE}{}{}
\patchcmd{\@setauthors}{\MakeUppercase}{\large}{}{}
\makeatother

\title{New congruences for 2-color partitions} 

\author[S. Chern]{Shane Chern}
\address{School of Mathematical Sciences, Zhejiang University, Hangzhou, 310027, China}
\email{\href{mailto:shanechern@zju.edu.cn}{shanechern@zju.edu.cn}; \href{mailto:chenxiaohang92@gmail.com}{chenxiaohang92@gmail.com}}

\date{}

\begin{document}

{\footnotesize\noindent \textit{J. Number Theory} \textbf{163} (2016), 474--481.\\
\doi{10.1016/j.jnt.2015.12.020}}

\bigskip \bigskip

\maketitle

\thispagestyle{empty}

\begin{abstract}
Let $p_k(n)$ denote the number of $2$-color partitions of $n$ where one of the colors appears only in parts that are multiples of $k$. We will prove a conjecture of Ahmed, Baruah, and Dastidar on congruences modulo $5$ for $p_k(n)$. Moreover, we will present some new congruences modulo $7$ for $p_4(n)$.


\MSC{Primary 11P83; Secondary 05A17.}
\end{abstract}

\section{Introduction}

A partition of a natural number $n$ is a nonincreasing sequence of positive integers whose sum equals $n$. For example, $\lambda=\{3,2,1\}$ is a partition of $6$ since $6=3+2+1$. Denote by $p(n)$ the number of partitions of $n$. It is well known that the generating function of $p(n)$ is given by
$$\sum_{n\ge 0}p(n)q^n=\frac{1}{(q;q)_\infty},\quad |q|<1,$$
where, as usual, we denote
$$(a;q)_\infty=\prod_{n\ge 0}(1-aq^n).$$

Let $p_k(n)$ be the number of $2$-color partitions of $n$ where one of the colors appears only in parts that are multiples of $k$. According to \cite{ABD2015}, its generating function is
\begin{equation}\label{eq:main}
\sum_{n\ge 0}p_k(n)q^n=\frac{1}{(q;q)_\infty(q^k;q^k)_\infty},\quad |q|<1.
\end{equation}
Furthermore, we can write $p_0(n):=p(n)$.

Many authors have studied the arithmetic properties of $p_k(n)$. For example, for $k=1$, Baruah and Sarmah \cite[Eq. (5.4)]{BS2013} showed that
$$p_1(25n+23)\equiv 0\pmod{5}.$$
For $k=2$, Chan \cite{Chan2010} obtained the following congruence modulo $3$
$$p_2(3n+2)\equiv 0\pmod{3},$$
while Chen and Lin \cite{CL2009} proved
$$p_2(25n+22)\equiv 0\pmod{5}$$
and
$$p_2(49n+t)\equiv 0\pmod{7}$$
where $t=15$, $29$, $36$, and $43$ by using the tool of modular forms. More recently, Ahmed, Baruah, and Dastidar \cite{ABD2015} found several new congruences modulo $5$:
$$p_k(25n+24-k)\equiv 0\pmod{5}$$
where $k=0$, $1$, $2$, $3$, $4$, $5$, $10$, $15$, and $20$. Meanwhile, they conjectured that the congruence also holds for $k=7$, $8$, and $17$.

In this paper, we shall give an affirmative answer to their conjecture, namely,
\begin{theorem}\label{th:1.1}
For any nonnegative integer $n$,
\begin{equation}\label{eq:th1.1}
p_k(25n+24-k)\equiv 0\pmod{5}
\end{equation}
where $k=7$, $8$, and $17$.
\end{theorem}

Moreover, we shall prove
\begin{theorem}\label{th:1.2}
For any nonnegative integer $n$,
\begin{equation}\label{eq:th1.2}
p_4(49n+t)\equiv 0\pmod{7}
\end{equation}
where $t=11$, $25$, $32$, and $39$.
\end{theorem}

Our method is based on a result of Radu and Sellers \cite{RS2011} relating to modular forms, which can be tracked back to \cite{Radu2009}.

\section{Preliminaries}

Let $\gamma=\begin{pmatrix}a & b\\c & d\end{pmatrix}\in SL_2(\mathbb{Z})=:\Gamma$. For a positive integer $N$, we define the congruence subgroup of level $N$ as
$$\Gamma_0(N)=\left\{\gamma\in SL_2(\mathbb{Z}): \gamma\equiv \begin{pmatrix}* & *\\0 & *\end{pmatrix}\pmod{N}\right\}$$
where ``$*$'' means ``unspecified.'' It is known that
$$[\Gamma:\Gamma_0(N)]=N\prod_{p\mid N}(1+p^{-1}).$$
Furthermore, we write
$$\Gamma_\infty=\left\{\left.\begin{pmatrix}1 & h \\0 & 1 \end{pmatrix}\ \right|\ h\in\mathbb{Z}\right\}.$$

Let $m,M,N$ be positive integers. We write $R(M)$ the set of integer sequences indexed by the positive divisors $\delta$ of $M$. For any $r\in R(M)$, it has the form $r=(r_{\delta_1},\ldots,r_{\delta_k})$ where $1=\delta_1<\cdots<\delta_k=M$ are positive divisors of $M$. Let $[s]_m$ be the set of all elements congruent to $s$ modulo $m$. Let $\mathbb{Z}_m^*$ denote the set of all invertible elements in $\mathbb{Z}_m$, and $\mathbb{S}_m$ denote the set of all squares in $\mathbb{Z}_m^*$. Define by $\overline{\odot}_r$ the map $\mathbb{S}_{24m}\times\{0,\ldots,m-1\}$$\to$$\{0,\ldots,m-1\}$ with
$$([s]_{24m},t)\mapsto [s]_{24m}\overline{\odot}_r t\equiv ts+\frac{s-1}{24}\sum_{\delta\mid M}\delta r_{\delta}\pmod{m}$$
where $t\in\{0,\ldots,m-1\}$. Now write $P_{m,r}(t)=\{[s]_{24m}\overline{\odot}_r t\ |\ [s]_{24m}\in\mathbb{S}_{24m}\}$. Let $\Delta^*$ be the set of tuples $(m,M,N,t,r=(r_\delta))$ which satisfy conditions given in \cite[p. 2255]{RS2011}\footnote{According to a private communication between the author and S. Radu, the last condition of $\Delta^*$ should read: ``for $(s,j)=\pi(M,(r_\delta))$ we have ( ($4\mid \kappa N$ and $8\mid Ns$) or ($2\mid s$ and $8\mid N(1-j)$) ) if $2\mid m$.'' In our cases, since $2$ does not divide $m=25$ or $49$, none of the two conditions need to be satisfied.}. Finally, for $\gamma=\begin{pmatrix}a & b \\c & d \end{pmatrix}$, $r\in R(M)$, and $r'\in R(N)$ we denote
$$p_{m,r}(\gamma)=\min_{\lambda\in\{0,\ldots,m-1\}}\frac{1}{24}\sum_{\delta\mid M}r_\delta\frac{\gcd^2(\delta(a+\kappa\lambda c),mc)}{\delta m}$$
and
$$p_{r'}^*(\gamma)=\frac{1}{24}\sum_{\delta\mid N} \frac{r'_\delta\gcd^2(\delta,c)}{\delta}$$
where $\kappa=\kappa(m)=\gcd(m^2-1,24)$.

Let
$$f_r(q):=\prod_{\delta\mid M}(q^{\delta};q^{\delta})_{\infty}^{r_\delta}=\sum_{n\ge 0}c_r(n)q^n$$
for some $r\in R(M)$. The following lemma (see \cite[Lemma 4.5]{Radu2009} or \cite[Lemma 2.4]{RS2011}) is a key to our proof.
\begin{lemma}\label{le:01}
Let $u$ be a positive integer, $(m,M,N,t,r=(r_\delta))\in\Delta^*$, $r'=(r'_\delta)\in R(N)$, $n$ be the number of double cosets in $\Gamma_0(N)\backslash\Gamma/\Gamma_\infty$ and $\{\gamma_1,\ldots,\gamma_n\}$ $\subset\Gamma$ be a complete set of representatives of the double coset $\Gamma_0(N)\backslash\Gamma/\Gamma_\infty$. Assume that $p_{m,r}(\gamma_i)+p_{r'}^*(\gamma_i)\ge 0$ for all $i=1,\ldots,n$. Let $t_{\min} := \min_{t'\in P_{m,r}(t)}t'$ and
$$v:=\frac{1}{24}\left(\left(\sum_{\delta\mid M}r_\delta+\sum_{\delta\mid N}r'_\delta\right)[\Gamma:\Gamma_0(N)]-\sum_{\delta\mid N}\delta r'_\delta\right)-\frac{1}{24m}\sum_{\delta\mid M}\delta r_\delta-\frac{t_{\min}}{m}.$$
Then if
$$\sum_{n=0}^{\lfloor v \rfloor}c_r(mn+t')q^n\equiv 0 \pmod{u}$$
for all $t'\in P_{m,r}(t)$, then
$$\sum_{n\ge 0}c_r(mn+t')q^n\equiv 0 \pmod{u}$$
for all $t'\in P_{m,r}(t)$.
\end{lemma}

It also readily follows by the binomial theorem that
\begin{lemma}[cf. {\cite[Lemma 1.2]{RS2011}}]\label{le:02}
Let $p$ be a prime and $\alpha$ a positive integer. Then
$$\frac{(q;q)_\infty^{p^\alpha}}{(q^p;q^p)_\infty^{p^{\alpha-1}}}\equiv 1\pmod{p^\alpha}.$$
\end{lemma}

\section{Proof of Theorem \ref{th:1.1}}

\subsection{The case $\boldsymbol{k=7}$}\label{subsec:th1.1.7}

Taking $k=7$ in \eqref{eq:main}, we have
\begin{equation}\label{eq:k7}
\sum_{n\ge 0}p_7(n)q^n=\frac{1}{(q;q)_\infty(q^7;q^7)_\infty}.
\end{equation}
It follows by Lemma \ref{le:02} that
\begin{equation}\label{eq:g7}
\sum_{n\ge 0}p_7(n)q^n\equiv\frac{(q;q)_\infty^4}{(q^5;q^5)_\infty(q^7;q^7)_\infty}=:\sum_{n\ge 0}g_{7,5}(n)q^n \pmod{5}.
\end{equation}

We first set
$$(m,M,N,t,r=(r_1,r_5,r_7,r_{35}))=(25,35,35,17,(4,-1,-1,0))\in\Delta^*.$$
By the definition of $P_{m,r}(t)$, we obtain
$$P_{m,r}(t)=\left\{t'\ |\ t'\equiv ts-(s-1)/3\ (\bmod\ m),0\le t'\le m-1,[s]_{24m}\in\mathbb{S}_{24m}\right\}.$$
One readily verifies $P_{m,r}(t)=\{17\}$. Now setting
$$r'=(r'_1,r'_5,r'_7,r'_{35})=(3, 0, 11, 0).$$
Let $\gamma_\delta=\begin{pmatrix}
1 & 0 \\
\delta & 1 
\end{pmatrix}$. By \cite[Lemma 2.6]{RS2011}, $\{\gamma_\delta:\delta\mid N\}$ contains a complete set of representatives of the double coset $\Gamma_0(N)\backslash\Gamma/\Gamma_\infty$. It is easy to verify that all these constants satisfy the assumption of Lemma \ref{le:01}. We thus obtain the upper bound $\lfloor v \rfloor=28$. Through \textit{Mathematica}, we verify that $g_{7,5}(25n+17)\equiv 0\pmod{5}$ holds for the first $29$ terms. It therefore follows by Lemma \ref{le:01} that
$$g_{7,5}(25n+17)\equiv 0\pmod{5}$$
holds for all $n\ge 0$. Now by \eqref{eq:g7} we have
$$p_7(25n+17)\equiv 0\pmod{5}$$
for all $n\ge 0$.

\subsection{The case $\boldsymbol{k=8}$}

Taking $k=8$ in \eqref{eq:main}, we have
\begin{equation}\label{eq:k8}
\sum_{n\ge 0}p_8(n)q^n=\frac{1}{(q;q)_\infty(q^8;q^8)_\infty}.
\end{equation}
By Lemma \ref{le:02}, one obtains
\begin{equation}\label{eq:g8}
\sum_{n\ge 0}p_8(n)q^n\equiv\frac{(q;q)_\infty^4}{(q^5;q^5)_\infty(q^8;q^8)_\infty}=:\sum_{n\ge 0}g_{8,5}(n)q^n \pmod{5}.
\end{equation}

In this case, we may set
\begin{align*}
(m,M,N,t,r=&(r_1, r_2, r_4, r_5, r_8, r_{10}, r_{20}, r_{40}))\\
&=(25,40,40,16,(4, 0, 0, -1, -1, 0, 0, 0))\in\Delta^*
\end{align*}
and
$$r'=(r'_1, r'_2, r'_4, r'_5, r'_8, r'_{10}, r'_{20}, r'_{40})=(0, 0, 0, 0, 14, 0, 0, 0).$$
We also obtain $P_{m,r}(t)=\{16\}$. One readily computes that $v$ is bounded by $\lfloor v \rfloor=42$. With the help of \textit{Mathematica}, we see that $g_{8,5}(25n+16)\equiv 0\pmod{5}$ holds up to the bound $\lfloor v \rfloor$. We conclude by Lemma \ref{le:01} and \eqref{eq:g8} that
$$p_8(25n+18)\equiv 0\pmod{5}$$
holds for all $n\ge 0$.

\subsection{The case $\boldsymbol{k=17}$}

Taking $k=17$ in \eqref{eq:main}, we have
\begin{equation}\label{eq:k17}
\sum_{n\ge 0}p_{17}(n)q^n=\frac{1}{(q;q)_\infty(q^{17};q^{17})_\infty}.
\end{equation}
Thanks to Lemma \ref{le:02}, one gets
\begin{equation}\label{eq:g17}
\sum_{n\ge 0}p_{17}(n)q^n\equiv\frac{(q;q)_\infty^4}{(q^5;q^5)_\infty(q^{17};q^{17})_\infty}=:\sum_{n\ge 0}g_{17,5}(n)q^n \pmod{5}.
\end{equation}

Here we set
$$(m,M,N,t,r=(r_1,r_5,r_{17},r_{85}))=(25,85,85,7,(4,-1,-1,0))\in\Delta^*$$
and
$$r'=(r'_1,r'_5,r'_{17},r'_{85})=(0, 0, 20, 0).$$
In this case $P_{m,r}(t)=\{7\}$. It is easy to get the upper bound $\lfloor v \rfloor=84$. Now we verify that $g_{17,5}(25n+7)\equiv 0\pmod{5}$ holds for the first $85$ terms by \textit{Mathematica}. It follows by Lemma \ref{le:01} and \eqref{eq:g17} that for all $n\ge 0$
$$p_{17}(25n+7)\equiv 0\pmod{5}.$$

\section{Proof of Theorem \ref{th:1.2}}

Taking $k=4$ in \eqref{eq:main}, we have
\begin{equation}\label{eq:k4}
\sum_{n\ge 0}p_{4}(n)q^n=\frac{1}{(q;q)_\infty(q^{4};q^{4})_\infty}.
\end{equation}
It follows by Lemma \ref{le:02} that
\begin{equation}\label{eq:g4}
\sum_{n\ge 0}p_{4}(n)q^n\equiv\frac{(q;q)_\infty^6}{(q^4;q^4)_\infty(q^{7};q^{7})_\infty}=:\sum_{n\ge 0}g_{4,7}(n)q^n \pmod{7}.
\end{equation}
To prove Theorem \ref{th:1.2}, it suffices to show
\begin{equation}\label{eq:7g4}
g_{4,7}(49n+t)\equiv 0\pmod{7}
\end{equation}
for $t=11$, $25$, $32$, and $39$.

We first prove the cases $t=11$, $25$, and $32$. Setting
\begin{align*}
(m,M,N,t,r=&(r_1, r_2, r_4, r_7, r_{14}, r_{28}))\\
&=(49,28,28,11,(6, 0, -1, -1, 0, 0))\in\Delta^*.
\end{align*}
We compute that $P_{m,r}(t)=\{11, 25, 32\}$. Now taking
$$r'=(r'_1, r'_2, r'_4, r'_7, r'_{14}, r'_{28})=(2, 0, 1, 0, 0, 0),$$
and choosing $\gamma$ as in Subsection \ref{subsec:th1.1.7}, we verify that all these constants satisfy the assumption of Lemma \ref{le:01}. We thus compute $\lfloor v \rfloor=13$. With the help of \textit{Mathematica}, we see that \eqref{eq:7g4} holds up to the bound $\lfloor v \rfloor$ with $t\in\{11, 25, 32\}$, and thus it holds for all $n\ge 0$ by Lemma \ref{le:01}.

Now we will turn to the case $t=39$. Again we set
\begin{align*}
(m,M,N,t,r=&(r_1, r_2, r_4, r_7, r_{14}, r_{28}))\\
&=(49,28,28,39,(6, 0, -1, -1, 0, 0))\in\Delta^*
\end{align*}
and
$$r'=(r'_1, r'_2, r'_4, r'_7, r'_{14}, r'_{28})=(1, 0, 1, 0, 0, 0).$$
In this case we have $P_{m,r}(t)=\{39\}$. One readily computes that $v$ is bounded by $\lfloor v \rfloor=11$. Similarly we verify the first $12$ terms of \eqref{eq:7g4} with $t=39$ through \textit{Mathematica}. It follows by Lemma \ref{le:01} that it holds for all $n\ge 0$.

This completes our proof of Theorem \ref{th:1.2}.

\subsection*{Acknowledgments}

The author would like to thank the anonymous referee for careful reading and useful comments.

\bibliographystyle{amsplain}

\end{document}